\documentclass[review]{elsarticle}
\usepackage[leqno]{amsmath}
\usepackage[all]{xy}
\usepackage[mathscr]{euscript}
\usepackage{amssymb}
\usepackage{amscd}
\usepackage{epsfig}
\usepackage{lineno,hyperref}
\modulolinenumbers[5]

\newcommand\ring[1]{\mathaccent23{#1}}
\def\Vr{\ring{V}} 

\newcommand{\pf}{\noindent{\em Proof: }}
\newcommand{\epf}{\hfill\hbox{\rule{3pt}{6pt}}\\}

\newcommand{\R}{{\mathbb R}}

\newcommand{\Hh}{\mathcal H}
\newcommand{\Pa}{\mathbf P}
\newcommand{\T}{\mathcal T}
\newcommand{\C}{\mathcal C}

\newtheorem{lemma}{Lemma}[section]
\newtheorem{theorem}[lemma]{Theorem}
\newtheorem{proposition}[lemma]{Proposition}
\newtheorem{corollary}[lemma]{Corollary}

\begin{document}

\begin{frontmatter}

\title{Combinatorial properties of triplet covers for binary trees}


\author[stefanaddress]{Stefan~Gr{\"u}newald}
\ead{stefan@picb.ac.cn}

\author[uaeaddress]{Katharina T.~Huber*}
\ead{K.Huber@uea.ac.uk}

\author[uaeaddress]{Vincent Moulton}
\ead{V.Moulton@uea.ac.uk}

\author[mikeaddress]{Mike Steel}
\ead{Mike.Steel@canterbury.ac.nz}

\address[stefanaddress]{CAS-MPG Partner Institute for Computational Biology, Chinese Academy of Sciences Key Laboratories of Computational Biology, Shanghai, China}
\address[uaeaddress]{School of Computing Sciences, University of East Anglia, Norwich, UK}
\address[mikeaddress]{Biomathematics Research Centre, University of Canterbury, Christchurch, NZ}

\begin{abstract}
It is a classical result that an unrooted tree $T$ having  positive real-valued edge lengths and no vertices of degree two can be  reconstructed from the induced distance between each pair of leaves.  Moreover, if each non-leaf vertex of $T$ has degree 3  then the number of distance values required is linear in the number of leaves.  A canonical candidate for such a set of pairs of leaves in $T$  is the following: for each non-leaf vertex $v$, choose a leaf in each of the three components of $T-v$, group these three leaves  into three pairs, and take the union of this set over all
choices of $v$. This forms a so-called  `triplet cover' for $T$. In the first part of this paper we answer an open question (from 2012) by showing that the induced leaf-to-leaf distances for any triplet cover for $T$ uniquely determine $T$ and its edge lengths. We then investigate the finer combinatorial properties of triplet covers.  In particular, we describe the structure of triplet covers that satisfy one or more of the following properties of being minimal, `sparse', and `shellable'. 
\end{abstract}


\end{frontmatter}

*Corresponding author

\section{Introduction}

Trees with a label  set $X$ of leaves play a central role in many areas of classification, such as systematic biology and linguistics.  In these settings, it is usually assumed that the non-leaf vertices of the tree have degree at least three, and that there is an assignment of a positive real-valued length to each edge of $T$.   A classical and important result from the 1960s and 1970s asserts that any such (unrooted)  tree $T$ with edge lengths is uniquely determined  from the induced leaf-to-leaf distances between each pair of elements of $X$.  This result is the basis of widely-used methods for inferring trees from distance data, such as the popular `Neighbor-Joining' algorithm \cite{sai}.

When the unrooted tree $T$ is binary (each non-leaf vertex has degree 3) then we do not require distance values for all of the $\binom{n}{2}$ pairs from leaf set $X$ (where $n=|X|\geq 3$),  since just $2n-3$ carefully selected pairs of leaves suffice to determine $T$ and its edge lengths (see \cite{GLM04};  more recent  results  appear in \cite{DHS12},  motivated by  the irregular distribution of genes across species in biological data).    This value of $2n- 3$ cannot be made any smaller, since a binary unrooted tree with $n$ leaves has $2n-3$ edges, and the inter-leaf distances are linear combinations of the corresponding $2n-3$ edge lengths (so, by linear algebra, these values cannot be uniquely determined by fewer than $2n-3$ equations).

There is a particularly natural way to select a subset of $\binom{X}{2}$ for $T$ when $T$ is binary.  Since each non-leaf
vertex  is incident with three subtrees of $T$, let us  (i) select a leaf from each subtree,  (ii) consider the three pairs of leaves
we can form from this triple, and  then (iii) take the union of these sets of pairs  over all non-leaf vertices of $T$. This process produces a `triplet cover' of $T$  (defined more precisely below), as illustrated in Fig.~\ref{fig1}.

A triplet cover need not be of this minimum size (i.e. of size $2n-3$) and in an earlier paper we characterized when it is \cite{HMS17}.  That paper also
established that in this case the resulting triplet cover is `shellable', complementing other recent work into
phylogenetic `lasso' sets \cite{DHS12,HS14}, as well as a Hall-type characterization of the median function on trees in \cite{DS09}. 

In this paper, we present three main new results. Our first  result (Theorem~\ref{firstmainthm}, which is a special case of Theorem~\ref{stefanthm})  answers in the affirmative a question that has been open since 2012, namely do the distances between leaves induced by a triplet cover on a binary tree with positive edge lengths determine the tree?  Our second main result (Theorem~\ref{secthm}) describes the structure of `sparse' triplet covers in terms of a 2-tree decomposition of a certain graph.  Our third main result (Theorem~\ref{patchthm}) provides a sufficient condition for a triplet cover to be `shellable'.
Along the way,  a number of other properties of triplet covers are derived.  We begin with some definitions.

\subsection{Definitions}

Let $X$ be a finite set $|X| \ge 3$. Given a set $\C$ of subsets of a set $Y$, we let $\bigcup \C = \bigcup_{t \in \C} t$, and we denote elements in $X \choose 2$ and $X \choose 3$ also by $xy$ and $xyz$, respectively, where $x,y,z \in X$.

Given a graph $G=(V,E)$, we let $V=V(G)$ denote its vertex set and $E=E(G)$ its edge-set.

\begin{figure}[htb]
\centering
\includegraphics[scale=0.8]{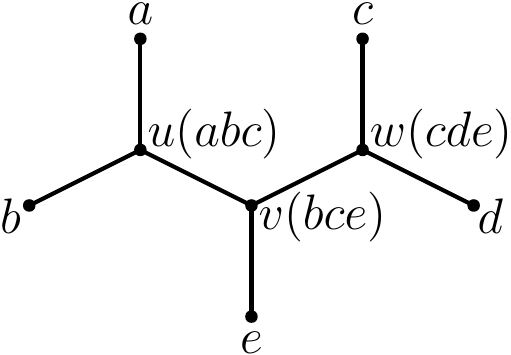}
\caption{A binary phylogenetic $X$--tree $T$ for $X=\{a,b,c,d,e\}$, and, in parentheses,  triples of elements from $X$ whose medians correspond to the indicated interior vertices.   The corresponding triplet cover is $\T=\{ab, ac, bc, be, ce, cd, de \}.$}
\label{fig1}
\end{figure}

\begin{itemize}

\item A {\em phylogenetic $X$--tree} is a tree $T=(V,E)$
that has leaf-set $X$ and for which the non-leaf vertices have degree at least 3.
If all non-leaf vertices of $T$ have degree exactly 3 then we say that $T$ is a {\em binary} phylogenetic $X$--tree
(or simply a binary phylogenetic tree when the leaf set is clear or not important). 

\item We let $\Vr  = \Vr(T) \subseteq V(T)$ denote the set of interior vertices of $T$.
If $T$ is binary then  $|\Vr|=|X|-2$. 

\item A {\em cherry} of a binary phylogenetic $X$--tree $T$ is a pair of leaves of $T$  that are adjacent to a common vertex.

\item Two phylogenetic $X$--trees $T=(V, E)$ and $T'=(V', E')$ are {\em isomorphic}, denoted $T \cong T'$,  precisely if there
is a graph isomorphism $\varphi$ from $T$ to $T'$ that 
sends leaf $x$ of $T$ to leaf $x$ of $T'$ for each $x \in X$. \\

\noindent Suppose that $T$ is a phylogenetic $X$-tree as above.

\item Given an element $x \in X$, where $|X| \ge 4$, we let $T-x$ denote the
phylogenetic $X$--tree which is obtained by
removing the leaf $x$ and the edge that contains it from $T$ 
and suppressing the remaining
degree 2 vertex.

\item Suppose that $\T$ is a subset of $X \choose 2$.
We say that a triple $abc \in {X \choose 3}$ {\em supports} a vertex
$v \in \Vr$ (relative to $\T$)
if we can select leaves $a,b,c \in X$, one from each
component of $T-v$, such that $ab, ac, bc \in \T$.

\item We call a subset $\T \subseteq {X \choose 2}$ a
{\em triplet cover} for $T$ if, for each element $v\in \Vr$, 
there is some element in ${X \choose 3}$ that supports $v$ 
(relative to $\T$).  An example is shown in Fig.~\ref{fig1}.  We call
each element in $\T$ a {\em cord}.

\item Given a subset $\T \subseteq {X \choose 2}$
and $x \in X$,
we let $\T^{-x} = \T - \{ xa \,: \, a \in X-\{x\} \mbox{ and } xa \in \T\}$.
In other words,  $\T^{-x}$ is the subset of $\T$ obtained by removing from $\T$
precisely those cords which contain $x$.

\end{itemize}

\section{Tree distances from any triplet cover determines the underlying binary tree}

The {\em triplet cover question} posed in \cite{DHS12} and discussed further in \cite{HS14},  asks the following.
Suppose that $\T$ is a triplet cover for a binary phylogenetic tree $X$--tree $T$ having strictly positive edge lengths.  Then  is $T$ the only phylogenetic $X$--tree (up to isomorphism) with strictly positive edge lengths that can induce the same distance values on pairs of leaves chosen from  $\T$?

It is fairly straightforward to show that the triplet cover question has an affirmative answer if we impose the additional restriction on each possible alternative tree $T'$  that $\T$ is also a triplet cover for $T'$ (Proposition 1 of \cite{DHS12}).   However, this additional restriction cannot be assumed {\em a priori} for the general question.
The following theorem provides an affirmative answer to the triplet cover question in general.

To state this theorem, we require one further definition. Given a phylogenetic $X$--tree $T$ and an assignment $\ell$ of strictly positive lengths to the edges of $T$,  let $d_{(T, \ell)}(x,y)$ denote the total length of the path between $x$ and $y$ in $T$ (i.e. the sum of the lengths of the edges in the unique path in $T$ that connects $x$ to $y$). 

\begin{theorem}
\label{firstmainthm}
Suppose that $\T$ is a triplet cover for a binary phylogenetic $X$--tree $T$, where $|X|\geq 3$, and that $\ell$ is an assignment of strictly positive edge lengths for $T$.
Then for any phylogenetic $X$--tree $T'$ and any assignment  $\ell'$ of strictly positive edge lengths for $T'$, if
$d_{(T, \ell)}(x,y) = d_{(T', \ell')}(x,y)$ for all $xy \in \T$ we have $T \cong T'$ and $\ell = \ell'$.
\end{theorem}

{\em Remarks:} It suffices to show that $T \cong T'$, since  $\ell = \ell'$ then follows easily (for example, by Proposition 1 of \cite{DHS12}).  Note also that
Theorem~\ref{firstmainthm} is not true when $|X|=2$, as we may take $T'$ to have a single vertex, and $T$ to have two vertices joined by an edge of arbitrary length (in this case $\T=\emptyset$).

We prove the other (main) part of the theorem by establishing a slightly more general result (Theorem~\ref{stefanthm}) which allows a more streamlined proof-by-contradiction argument based on the assumption that a (minimal) counterexample exists.

The proof of the following result (Proposition~\ref{construct}) does not  require or use Theorem~\ref{firstmainthm} (nor is that theorem implied by Proposition~\ref{construct}). However, the two results are complementary since 
Theorem~\ref{firstmainthm} ensures that the reconstructed tree (namely $T$) produced by the algorithm described in  Proposition~\ref{construct}  is the only tree that can realize $d$ for all elements $xy \in \T$.

  \begin{proposition}
  \label{construct}
  There is a polynomial-time algorithm that will reconstruct any binary phylogenetic $X$--tree $T$ together with a strictly positive edge length assignment $\ell$ given the values of
  $d=d_{(T, \ell)}$ on elements $xy \in \T$ where $\T$ is a triplet cover for $T$.
  \end{proposition}

\pf
 To this end, for any leaf $x$ of $X$ let  $$\lambda(x) =\frac{1}{2} \min\{d(x,z) + d(x, z') -d(z,z'): xz, xz' \in \T, z \neq z'\}.$$
  Notice that the set in the definition of $\lambda(x)$ is  non-empty since $\T$ is a triplet cover of $T$ (consider the vertex of $T$ adjacent to leaf $x$) and so $\lambda(x)$ can
be determined (in polynomial-time) from the values of $d$ on $\T$. Moreover, $\lambda(x)$ is the length of the pendant edge of $T$ incident with leaf $x$.

CLAIM:  $x$ and $y$ form a cherry of $T$ if and only if $xy \in \T$ and $d(x,y) = \lambda(x)+\lambda(y)$.

\noindent {\em Proof of claim:} First suppose that $x$ and $y$ form a cherry of $T$. We have $xy \in \T$ since $\T$ is a triplet cover for $T$ (consider the  interior vertex of $T$ that is adjacent to $x$ and $y$).  Moreover, $d(x,y) = \lambda(x)+\lambda(y)$.  Conversely if $x$ and $y$ do not form a cherry of $T$ and if $xy \in \T$ then $d(x,y) > \lambda(x)+\lambda(y)$ since there is at least one interior edge, having strictly positive length, in the path between $x$ and $y$.  This establishes the claim.

It follows form the claim that from $\T$ and $d|_\T$ we can identify a cherry $x,y$ of $T$, as well as the length of the pendant edges incident with $x$ and $y$ (namely, $\lambda(x)$ and $\lambda(y)$).
Let us now remove  $x$ from $X$ to give a reduced label set $X'=X-x$, remove $xy$ from $\T$ and  replace each pair $xz \in \T$ by a pair $yz$ to  obtain a modified set $\T' \subseteq \binom{X'}{2}$. For any pair $xz \in \T$ that has been replaced by $yz \in \T'$ set $d'(y,z) = d(x,z)+\lambda(y)-\lambda(x)$, and let $d'$ coincide with $d$ for all other elements of $\T'$.
  Then $\T'$ is now a triplet cover for the binary tree $T-x$ with its induced (strictly positive) edge lengths.
By induction on $|X|$ we can continue this cherry identification and deletion process until we obtain a binary phylogenetic tree on just three leaves. By reversing this process the original tree $T$ and its
edge lengths $\ell$ can then be reconstructed (in polynomial time).
\epf

\subsection{Proof of Theorem~\ref{firstmainthm}}

As mentioned, we need only prove that a binary phylogenetic tree $T$ is uniquely determined (up to isomorphism) by the values of $d(x,y)=d_{(T, l)}(x,y)$ (for some strictly positive assignment $l$ of edge lengths to $T$)  for pairs $x,y$ in a triplet cover of $T$ (the uniqueness of the edge lengths was shown in \cite{DHS12} and is straightforward).  Also, it is convenient for the proof to work with a slightly more general class of trees than binary phylogenetic trees, namely `binary $X$--trees'.  To define this class, recall that an  {\em $X$--tree} $T$ is a tree for which each vertex of degree at most 2 is labelled by at least one element of $X$,
and each element of $X$ labels exactly one vertex of the tree. 
Thus $T$ may have unlabelled vertices, but these must have degree 3 or more, and $T$ may have vertices labelled by more than one element of $X$. We refer to $X$ as the {\em label set} of $T$. 
 Thus, a phylogenetic $X$--tree is an $X$--tree in which only the leaves are labelled, and each leaf is labelled by only a single element of $X$. 

Deleting an edge $e$ of a $X$--tree and  considering the two connected components of $T-e$ partitions $X$ into two nonempty sets. Any such bipartition of $X$ consisting of two blocks, say $A$ and $B$, is called an {\em $X$--split} and we denote it by writing $A|B$.   Moreover, the set of splits arising from an $X$--tree $T$ in this way (by cutting edges of $T$) is {\em pairwise compatible}, which means that, for any two splits  $A|B$ and $A'|B'$ in the set, at least one of the four intersections $A\cap A', A \cap B', B \cap A', B \cap B'$ is empty.  It is a classical result (due to Peter Buneman) that  $X$--trees (up to a natural notion of isomorphism) are in bijective correspondence with the sets of pairwise compatible $X$--splits, with each $X$--split in the set corresponding to a unique edge of the associated tree.   Thus two $X$--trees $T$ and $T'$ are  isomorphic,  written $T \cong T'$,  precisely if they have the same set of $X$--splits. This definition agrees with the earlier notion of isomorphism when restricted to binary phylogenetic $X$--trees. 

\begin{figure}[htb]
\centering
\includegraphics[scale=1.0]{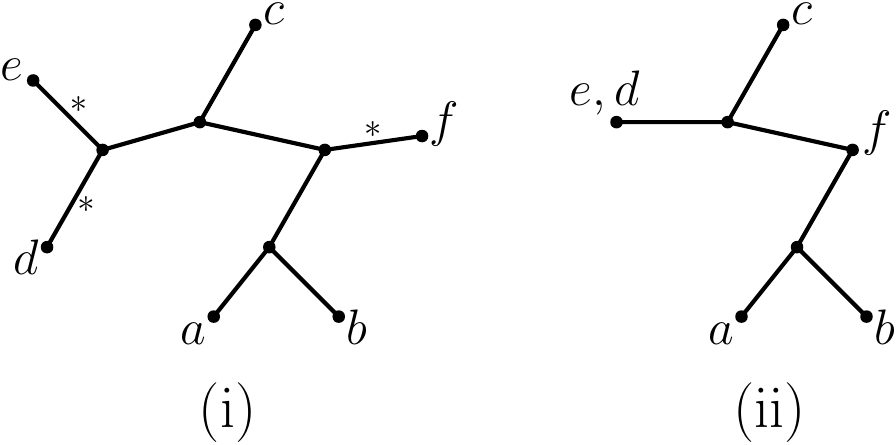}
\caption{(i) A binary phylogenetic $X$--tree for $X=\{a,b,c,d,e,f\}$. If we collapse the three edges indicated by * we obtain the binary $X$--tree shown in (ii). }
\label{binfig}
\end{figure}

We will say that an $X$--tree $T$ is  {\em binary} if  (i) each vertex of degree at most 1 is labelled by one or two elements of $X$, (ii) each vertex of degree 2
is labelled by exactly one element of $X$, and (iii) all remaining vertices are unlabelled and have degree 3.
Note that if we take any binary phylogenetic $X$--tree and collapse any (possibly empty) subset of pendant edges of this tree, we obtain a binary $X$--tree; conversely, each
binary $X$--tree $T$ is obtained from a unique binary phylogenetic $X$--tree $T_B$ by collapsing a unique subset of pendant edges of $T$.   Note that $T_B$ is the $X$--tree
whose set of $X$--splits consists of the $X$--splits of $T$ together with any {\em trivial $X$--splits} (i.e. $x|X-x, x \in X$) not already present in $T$.
This is illustrated in Fig.~\ref{binfig}.

It is straightforward to extend the notion of triplet cover from a binary phylogenetic $X$--tree to a binary $X$--tree: We say that a subset $\T$ of $\binom{X}{2}$ is a {\em triplet cover} of a binary $X$--tree $T$,
if $\T$ is a triplet cover (in the usual sense) of the associated binary phylogenetic $X$--tree $T_B$.

Given an $X$--tree $T$ and a function $\ell$ that assigns strictly positive lengths to each edge of $T$, let $d_{(T, \ell)}: X \times X \rightarrow \R^{\geq 0}$ be the induced
distance function on $X$ in which $d_{(T, \ell)}(x,y)$ is the sum of the lengths of the edges on the
(unique) path in $T$ connecting the vertices of $T$ labelled by $x$ and $y$.
Notice that $d_{(T, \ell)}$ takes the value 0 if and only if $x$ and $y$ label the same vertex of $T$. The function $d_{(T, \ell)}$ is non-negative, symmetric and satisfies the
triangle inequality, and so $d_{(T, \ell)}$ is a pseudometric on $X$.

\begin{theorem}
\label{stefanthm}
Suppose that $\T$ is a triplet cover for a binary $X$--tree $T$, where $|X|\geq 3$, and that $\ell$ is an assignment of strictly positive edge lengths for $T$.
Then for any $X$--tree $T'$ and any assignment $\ell'$ of  strictly positive edge lengths for $T'$, if
$d_{(T, \ell)}(x,y) = d_{(T', \ell')}(x,y)$ for all $xy \in \T$ we have $T \cong T'$.
\end{theorem}

\pf
The theorem is readily verified for $|X| = 3$. Suppose there is a counterexample to Theorem~\ref{stefanthm} when $|X|\geq 4$.
In that case, we can select a counterexample --- say $\T, T, T'$ --- for which
(i) $|X|\geq 4$ is minimal (call this minimal value $n$) and (ii) within all counterexamples with $|X|=n$  the sum of the number of edges in $T$ and in $T'$ is minimal.
We will show that we can then always construct another counterexample that either (i) has a label set of size $n-1\geq 3$, or has a smaller total number of edges across the two trees,  contradicting the minimality assumptions, or (ii) has a label set of size 3, for which the result holds.
Such a contradiction implies that no counterexample $\T, T, T'$ to Theorem~\ref{stefanthm} can exist.


We first establish the following claims:
\begin{itemize}
\item[{\bf (i)}] There is no trivial $X$--split present in both  $T$ and $T'$.
\item[{\bf (ii)}] $T'$ does not contain any trivial $X$--split.
\end{itemize}
To establish Claim (i), suppose that some trivial $X$--split, say $x|X-x$ (for some $x \in X$)  is present in both $T$ and $T'$.
Let $e_x$ and $e'_x$ denote the pendant edges of $T$ and $T'$ that are incident with $x$,  let $\ell(e_x), \ell'(e'_x)$ denote their lengths, and let $\ell_x = \min\{ \ell(e_x),  \ell'(e'_x) \}$.
If $\ell(e_x) = \ell'(e'_x)$ then collapse the pendant edges $e_x$ and $e'_x$, while if $\ell(e_x) \neq \ell'(e'_x)$ then collapse the shorter of the two pendant edges
$e_x$ and $e'_x$, and reduce the length of the other pendant edge by $\ell_x$.
The resulting pair of modified trees still consists of a binary $X$--tree and an $X$--tree, both with strictly positive edge lengths, and the distance between $x$ and any other element  of $X$ in both $T$ and in $T'$ has been reduced by $\ell_x$ while all other distances remain the same. Thus,  the sum of the number of edges in $T$ and $T'$ has been reduced by at least 1.  However the two modified trees, along with the original triplet cover $\T$ provide a smaller counterexample (in terms of the sum of the number of edges in $T$ and $T'$)  violating the minimality assumption. This establishes Claim (i).

To establish Claim (ii), assume that there is a trivial $X$--split $x|X-x$ (for some $x \in X$) present in $T'$.  By Claim (i), $x|X-x$ is not present in $T$, and so either $x$ labels an interior vertex of $T$, or else $x$ together with another element of $X$, say $y$, labels a leaf of $T$. In either case, since $\T$ is a triplet cover for $T$ there is a triplet $wxy$ with $\{wx, wy, xy\} \subseteq \T$ and with
 $d_{(T, \ell)}(w, x)+d_{(T, \ell)}(x,y)=d_{(T, \ell)}(w,y).$  But then (since $d_{(T, \ell)}(x',x'') = d_{(T', \ell')}(x',x'')$ for all $x'x'' \in \T$) it follows that $T'$ does not have $x$ as the sole label of one of its leaves (i.e. $T'$ does not contain $x|X-x$). Hence, $T'$ does not contain any trivial $X$--split.  This establishes Claim (ii).

Now, let $v$ be any leaf of $T'$,  and let $Y \subset X$ denote its label set
(note that $T'$ cannot be a single-vertex tree since otherwise $T$ would be also, and so $|X|\leq 2$, which violates our assumption that $n>3$).
By Claim (ii), $|Y|>1$.  Thus one of the following two cases must apply:
 \begin{itemize}
\item[]  {\bf Case A:} $T$ contains all trivial $X$--splits $y|X-y$ for all $y \in Y$.
\item[] {\bf Case B:} There exists some $y \in Y$ for which the trivial $X$--split $y|X-y$ is absent from $T$.
\end{itemize}
(We will show that neither case can arise, which will furnish the required contradiction).

 In Case A, let $e'_v$ denote the pendant edge of $T'$ that is incident with $v$, and for each $y \in Y$, let $\ell(e_y)$ be the length of the pendant edge of $T$ incident with the leaf labelled by (only) $y$. Let
 $$\ell_Y = \min[ \{\ell(e_y): y \in Y\} \cup \{\ell'(e'_v\}].$$
 If $\ell'(e'_v) = \ell_Y$ then collapse edge $e'_y$ in $T'$, otherwise reduce the length of $e'_v$ by $\ell_Y$.
 Similarly, for  each $y\in Y$ for which $\ell(e_y) = \ell_Y$ collapse edge $e_y$ of $T$, otherwise
 reduce the length of $e_y$ by $\ell_Y$.
Again the resulting modified  pair of trees still consists of a binary $X$--tree and an $X$--tree, both with strictly positive edge lengths, but either $T$ or $T'$ has at least one (pendant) edge fewer than before. Moreover, the distance between two distinct elements $x, x'$ of $X$  in either tree is either unchanged (if neither $x$ nor $x'$ is in $Y$) or is reduced by $\ell_Y$ for both trees, when $|\{x, x'\} \cap Y|=1$.  For the remaining case where $x,x' \in Y$ the distances
in the modified trees may differ, however we also have that $xx' \not\in \T$, since the distance between $x$ and $x'$ in $T$ and in $T'$ is different (it is zero in $T'$ and non-zero in $T$). It follows that the modified trees again provide a smaller counterexample (in terms of the sum of the number of edges in $T$ and in $T'$)  violating the minimality assumption. This shows that Case A cannot arise.

For Case B, there exists some element $y \in Y$ that labels either an interior vertex of $T$, or else $y$ together with another element of $X$, say $x$, labels a leaf of $T$. In either case, since $\T$ is a triplet cover for $T$ there is a triplet $xyz $ with $\{xy, xz, yz\} \subseteq \T$ and with
 $d_{(T, \ell)}(x,y )+d_{(T, \ell)}(y,z)=d_{(T, \ell)}(x,z).$
Thus, since $v$ is a leaf of $T'$ we must have that $x$ or $z$ is in $Y$. Without loss of generality we
may assume that $x  \in Y$. Then since $xy \in \T$ it follows that $d_{(T, \ell)}(x,y) = d_{(T', \ell')}(x,y) =0$.
In particular, $\{x,y\}$ is contained in the label set of a leaf $u$ of $T$.

Finally, delete  label $x$ to obtain a set $X' = X-\{x\}$ of size $|X|-1$, and
 form modified $X'$--trees $\tilde{T}$ and $\tilde{T'}$  from $T$ and $T'$ (respectively) by deleting the label $x$ from $u$ and $v$.
 Then $\tilde{T}$ is a binary $X'$--tree,  $\tilde{T'}$ an $X'$--tree, and both trees inherit strictly positive edge lengths from $T$ and $T'$.
Consider the modified set $\tilde{\T} \subseteq \binom{X'}{2}$
obtained from $\T$ by deleting $xy \in \T$ and replacing each remaining occurrence of $x$ by $y$ in any pair $xz \in \T$. Then $\tilde{\T}$ is a triplet cover for $\tilde{T}$,
and if $x'x'' \in \tilde{\T}$ then the distance between $x'$ and $x''$ in $\tilde{T}$ is the same as it is in $\tilde{T'}$. Consequently, $(\tilde{\T}, \tilde{T}, \tilde{T'})$ provides
 a counterexample to Theorem~\ref{stefanthm}.
However, this new counterexample has a label set $X'$ of size $|X'|=n-1$ which is one less than the starting counterexample $(\T, T, T')$.  If $n-1=3$ this is impossible, since the Theorem holds when the label set has size 3, while if $n-1\geq 4$ we have violated the minimality assumption regarding the choice of $(\T, T, T')$. This shows that Case B cannot arise, thereby completing the proof.
\epf

\section{Properties of minimal triplet covers}

In this section, we collect together a  number of definitions, observations, and results (extending earlier work from \cite{HMS17}) that are required for establishing some further main results later in the paper.

\subsection{Preliminaries}

In the remainder of this paper, unless stated to the contrary, we will assume that $T$ refers to a {\bf binary} phylogenetic $X$--tree.  
We will also write $\Vr$ for $\Vr(T)$ (the set of interior vertices of $T$) when $T$ is clear.

Suppose that $\T$ is a triplet cover of $T$. The following terminology and result is from  \cite{HMS17}.
 \begin{itemize}
 \item
 For $x \in X$ the {\em multiplicity $\mu_{\T}(x)$ of $x$ (relative to $\T$)}
is the number of cords in $\T$
that contain $x$. The {\em multiplicity of $\T$} is $\mu(\T) = \min_{x \in X} \mu_{\T}(x)$.
\item 
$|\T| \ge 2|X|-3$ (for a direct proof see \cite[Proposition 3]{HMS17}).
We call $\T$ {\em minimum} if  $|\T| = 2|X|-3$.

\item We call a triplet cover $\T$ of $T$ {\em minimal} if $\T-\{t\}$ is not
a triplet cover for $T$ for all $t \in \T$.

\end{itemize}

Note that there exist triplet covers that are minimal but not minimum (an example is given in  Fig.~\ref{minimal}).

The following lemma summarizes some results 
established in \cite{HMS17}  (namely, Corollary 1, Proposition 2, and Corollary 2 of that paper, respectively).

\begin{lemma}
\label{lemo2}
Suppose that $\T$ is a minimal triplet cover of $T$. Then

\begin{itemize}
\item[{\rm (i)}] $2|X|-3 \le |\T| \le 3|X|-6.$

\item[{\rm (ii)}]  $2 \le \mu(\T) \le 5$.

\item[{\rm (iii)}] If $\T$ is a minimum triplet cover, then $\mu(\T)=2$.
\end{itemize}
\end{lemma}

\bigskip

Given a subset $\T \subseteq {X \choose 2}$
and $v \in \Vr$,
we let $S_v(\T)$ be the subset of ${X \choose 3}$ which
consists of precisely those triples which support $v$.
We call $S_v(\T)$ the {\em support of $v$ (relative to $\T$)}.

Note that if $abc \in S_v(\T)$, some $v \in \Vr$,
then $v = {\rm med}_T(a,b,c)$, where for all $xyz \in {X \choose 3}$, 
${\rm med}_T(x,y,z)$ denotes
the {\em median} of $x,y,z$ (i.e. the unique vertex 
that lies on all of the shortest paths between $x$, $y$ and $z$).

\begin{lemma}
\label{lemo}
\mbox{}
\begin{itemize}
\item[{\rm (i)}]   $\T$ is a triplet  cover of $T$ if and only if $|S_v(\T)| \ge 1$ for all $v \in \Vr$.
\item[{\rm (ii)}]  If $\T$ is a triplet cover of $T$
and $v, w \in \Vr$ distinct, then $S_v(\T) \cap S_w(\T) = \emptyset$.
\end{itemize}
\end{lemma}

\pf
The proof of Part (i) is straightforward. For Part (ii),  if this were not the case, then
for $xyz \in S_v(\T) \cap S_w(\T)$ we would have $v={\rm med}_T(x,y,z)=w$,
a contradiction.
\epf

\subsection{The cover graph $\Gamma(\T)$ and triplet set $\mathcal C(\T)$}

Given subset $\T \subseteq {X \choose 2}$, 
the {\em cover graph (of $\T$)}, denoted $\Gamma(\T)$, is the graph with vertex set $X$ and
edge set $\T$.  This graph,  introduced in \cite{DHS12},  has played an important role in subsequent papers \cite{DHS14}, \cite{HS14}, \cite{HMS17}.

We now consider a set which will be useful for understanding the triangles (i.e. 3--cycles)  in a
cover graph. Given a subset $\T \subseteq {X \choose 2}$, we define 
$$
\C(\T) =  \dot\bigcup_{v \in \Vr} S_v(\T) \subseteq {X \choose 3}.
$$
Note that the union in this definition is disjoint by Lemma~\ref{lemo}(ii). In addition
we note some other useful properties of the set $\C(\T)$.

\begin{lemma}
	\label{lemo3}
	\mbox{}
	\begin{itemize}
		\item[{\rm (i)}] If  $\T$ is a triplet cover for $T$, then $\bigcup \C(\T) =X$, and $|\C(\T)| \ge |X|-2$.
		\item[{\rm (ii)}] If $\T$ is a minimal triplet cover, then every cord in $\T$ is a subset of
		some element of $\C(\T)$.
	\end{itemize}
\end{lemma}

\pf
For Part (i),  if  $x \in X$, then let $v$ be the vertex in $T$ adjacent to $x$. If
$A \in S_v(\T)$, then clearly $x \in A$. The inequality now follows from Lemma~\ref{lemo}(ii).
For Part (ii), suppose that $\T$ is minimal and 
that there is a cord $xy \in \T$ that is not a subset of any element
in $ \C(\T)$.
Then, for all $v \in \Vr$, $xy \in \T$ is not a subset of any element of $S_v(\T)$. It
follows by Lemma~\ref{lemo}(i) that $\T -\{xy\}$ is a triplet cover for $T$, a contradiction.

\epf

We now collect together some important properties of the cover graph.

\begin{theorem}
\label{covpro}
Suppose that $\T$ is a triplet cover for $T$.

\begin{itemize}
\item[{\rm (i)}]  The triangles in the cover graph $\Gamma(\T)$
are in bijective correspondence with the elements of $\C(\T)$.

\item[{\rm (ii)}] $\Gamma(\T)$ is 2-connected.\footnote{The connectivity of
	$\Gamma(\T)$ also follows from \cite{DHS12} (Proposition 1 and Corollary 3).}

\item[{\rm (iii)}]  If  $\T$ is a minimal triplet cover for $T$, then
every cord in $\T$ is the edge of some triangle in $\Gamma(\T)$.\\

\end{itemize}
\end{theorem}

\pf
{\em Part (i):}
Suppose $xyz \in \C(\T)$, and so $xyz \in S_v(\T)$, for some $v \in \Vr$. Then clearly $x,y,z$ is a triangle in
$\Gamma(\T)$, since $xy, yz, xz \in \T$.  Thus, we have a map $\psi$ that takes elements in $\C(\T)$  to triangles in $\Gamma(\T)$.
Clearly this map is injective. Moreover, if $x,y,z$ is a triangle in $\Gamma(\T)$, then
for $v = {\rm med}_T(x,y,z)$, we have $xyz \in S_v(\T)$ and $\psi(xyz)= x,y,z$.
Thus $\psi$ is surjective.

\bigskip

\noindent{\em Part (ii):} The statement clearly holds if $|X| =3$, and so 
we assume $|X|\ge 4$.
Suppose $x \in X$. Let $v \in V(T)$ be the vertex in $T$ adjacent
to $x$. Let $w \neq x$ be a vertex adjacent to $v$ in $T$, and let
$T_w$ be the tree which is the connected component of $T$
minus the edge $\{v,w\}$ that contains $w$. Let $Y_w = V(T_w) \cap X$.
We claim that the graph induced by $\Gamma(\T)$ on $Y_w$ is connected.

Consider $T_w$ as being a rooted, directed tree, with root $w$
and all edges directed away from $w$. For $u \in V(T_w)$
we let $Y_u$ denote the set of leaves $x$ in $X$ for which $u$ lies on the path in $T_w$ from the root vertex $w$  to $x$.
We now prove the claim using induction on $|Y_u|$. If $|Y_u| =1$, then
clearly the graph induced  by $\Gamma(\T)$ on $Y_u$ is connected.
Now, suppose $u \in V(T_w)$ with $|Y_u| > 1$.
Let $u', u''$ be the children of $u$ in $T_w$.
By induction, we can assume that the graphs induced by
$\Gamma(\T)$ on $Y_{u'}$ and on $Y_{u''}$ are connected.
But as $\T$ is a triplet cover for $T$, there must
exist a cord $y' y'' \in \T$ with $y' \in Y_{u'}$ and $y'' \in Y_{u''}$.
So the graph induced by $\Gamma(\T)$ on $Y_{u} = Y_{u'} \cup Y_{u''}$ is connected.
The claim now follows as $w \in V(T_w)$.

Now, let $xyz$ be a triple in the support $S_v(\T)$ of $v$, for some $y,z\in X$,
which must exist as $\T$ is a triplet cover for $T$.
Let $w',w''$ denote the two vertices in $T$ adjacent to $v$ that
are not equal to $x$. Then as $xyz  \in S_v(\T)$,
it follows without loss of generality that $y \in Y_{w'}$ and $z \in Y_{w''}$. Moreover,
as the graphs induced by $\Gamma(\T)$ on
$Y_{w'}$ and $Y_{w''}$ are both connected by the above claim, $x \not\in Y_{w'}, Y_{w''}$
and $\{y,z\}$ is an edge in $\Gamma(\T)$,
it immediately follows that the graph $\Gamma(\T)-x$ obtained
by removing $x$ and the edges which contain it 
is connected.  Note that such edges must exist as 
$xy,xz\in\T$.
As the selection of $x \in X$ was arbitrary, if follows that
$\Gamma(\T)-x$ is connected for all $x \in X$, i.e. $\Gamma(\T)$ is 2-connected.

\bigskip

\noindent {\em Part (iii):} This follows by Part (i) and Lemma~\ref{lemo3}(ii).

\epf

\subsection{Sparse triplet covers and Hall-type subsets of triples}

\begin{figure}[ht]
    \begin{center}
        \includegraphics[width=8cm]{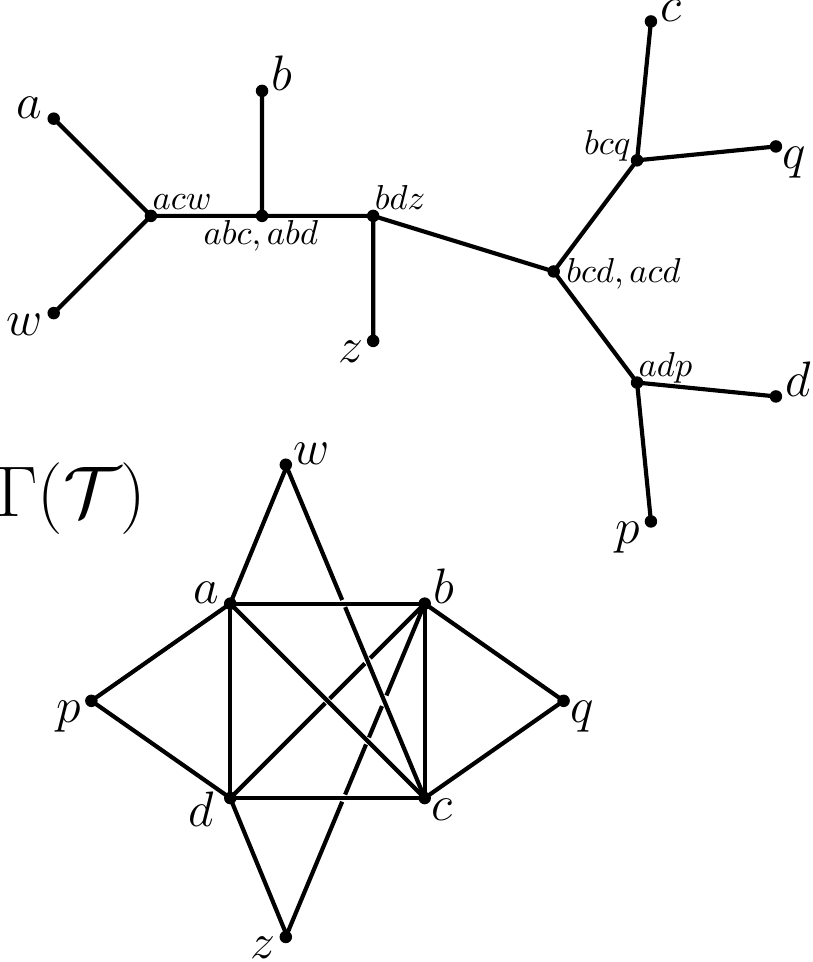}
    \end{center}
    \caption{Top:  A minimal triplet cover that is not sparse. Bottom: The associated cover graph for $\T$.}\label{minimal}
\end{figure}

We call a triplet cover $\T$ of $T$ {\em sparse} if  $|\C(\T)| = |X|-2$ (cf. Lemma~\ref{lemo3}(i)).
Note that if $\T$ is a sparse cover of $T$, then the function $f_{\T}: \C(\T) \to \Vr$,
which maps $xyz \in \C(\T)$ to ${\rm med}_T(x,y,z)$ is a bijection.
From this observation, it is possible to construct minimal triplet covers that are not sparse (an example is provided in Fig.~\ref{minimal}).

Note also that there are sparse triplet covers that are not minimal (an example is to
add a new cord $eh$ into the triplet cover $\T$ in Fig.~\ref{sparse}).

We say that a subset $\C \subseteq {X \choose 3}$ is of {\em Hall-type} if $\bigcup \C = X$ and $\C$ satisfies the
following property ({\em cf.} \cite{DS09}):  For all non-empty subsets $\C' \subseteq \C$,
$$
\left|\bigcup \C'\right| \ge |\C'| + 2.
$$
For example, for the triplet cover $\T$ in Fig.~\ref{sparse}, the set $\C(\T)$ is of Hall-type.

\begin{lemma}
\label{lemo4}
Suppose $\T$ is a triplet cover of $T$. If $\T$ is sparse, then $\C(\T)$ is
of Hall-type.
\end{lemma}

\pf If $\T$ is sparse then the map $f_{\T}$ is a bijection. Hence,
since $\bigcup \C(\T)=X$ by Lemma~\ref{lemo3}(ii),
it follows from \cite[Theorem 1.1]{DS09} that $\C(\T)$ is of Hall-type.
\epf


\subsection{Sections}

A subset $\C \subseteq \C(\T)$
is called a {\em section (of $\C(\T)$)} if $|\C \cap S_v(\T)| = 1$ for all $v \in \Vr$.
Note that $\C(\T)$ always contains a section and that 
if $\C$ is a section, then $|\C| = |X|-2$ 

We now define a set whose properties will be useful later on.
Given a subset $\C \subseteq {X \choose 3}$, we define the {\em cord set of $\C$} to be
$$
Co(\C) = \{ xy \in {X \choose 2} \,:\, xy \subset A \mbox{ some } A \in \C\}.
$$
Note that clearly $Co(\C(\T)) \subseteq \T$, but that in general $\T$ is not necessarily
a subset of $Co(\C(\T))$ (for example, if we
add in a new cord $eh$ to the triplet cover $\T$ in Fig.~\ref{sparse}, we obtain
a new triplet cover $\T'$ for which $\T'$ is not a subset of $Co(\C(\T'))$).
Moreover, if $\C$ is a section of $\C(\T)$, then $Co(\C) \subseteq \T$.

\begin{proposition}
\label{lemtrip}
Suppose that $\T$ is a triplet cover for $T$.  The following hold:
\begin{itemize}
\item[{\rm (a)}] If $\T$ is a minimal triplet cover for $T$, then $Co(\C(\T)) = \T$.
\item[{\rm (b)}]  $\T$ is sparse if and only if $\C(\T)$ is of Hall-type.
\item [{\rm (c)}]
\begin{itemize}
\item[{\rm (i)}] If $\C$ is a section of $\C(\T)$, then $Co(\C)$ is a triplet cover for $T$. In particular, $\bigcup\C = X$.
\item[{\rm (ii)}] If $\C$ is a section of $\C(\T)$, then $\C$ is of Hall-type.
\item[{\rm (iii)}] $\T$  is minimal if and only if $Co(\C) = \T$ for every section $\C$ of $\C(\T)$.     
\item[{\rm (iv})]  $\T$ is sparse if and only if $\C(\T)$ has a unique section.
\end{itemize}

\end{itemize}

\end{proposition}

\pf
{\em Part (a):}  Since $\T$ is minimal, for all $xy \in \T$, $xy \subset A$ for
all $A \in S_v(\T)$, some $v \in \Vr$ (otherwise we could remove $xy$ from $\T$
and still have a triplet cover).
So $\T \subseteq Co(\C(\T))$. As remarked above, the reverse inequality is obvious.

\bigskip

\noindent{\em Part (b):}
By Lemma~\ref{lemo4} it suffices to prove that if $\C(\T)$ is of Hall-type, then $\T$ is sparse.
Suppose for contradiction that $\T$ is not sparse. Let $\C$ be
a section of $\C(\T)$, so that $|\C|=|X|-2$ and so (since   $\C(\T)$ is of Hall-type) we have $\bigcup \C =X$. 
Since $\T$ is not sparse $\C(\T) > |X| -2$. Hence there is some $t \in \C(\T)$
that is not in $\C$. Let $\C' = \C \cup \{t\}$. Since $\C(\T)$ is of Hall-type and $\C' \subseteq \C(\T)$, 
$$
|X|= |\bigcup \C'| \ge |\C'|+2 = (|X|-2+1)+2 = |X|+1,
$$
a contradiction.

\bigskip

\noindent{\em Part (c-i):}  If $v \in \Vr(T)$, then there exists 
$xyz \in \C$, $x,y,z\in X$, with $xyz \in S_v(\T)$,
and so $xy, yz, zx \in Co(\C)$. Hence $S_v(Co(\C)) \neq \emptyset$. 
The statement now follows from Lemma~\ref{lemo}(i).\\

{\em Part (c-ii):} 
We can think of a section $\C$ of $\C(\T)$ as being a bijective map $f_{\C}: \Vr \to \C$
which for each $v \in \Vr$ selects some element in $S_v(\T)$ (the
inverse of $f_{\C}$ is the map which takes each $xyz \in \C$,
$x,y,z\in X$,  to ${\rm med}_T(x,y,z)$).
Statement (ii) now follows immediately from \cite[Theorem 1.1]{DS09}.\\

{\em Part (c-iii):}  Suppose $\T$ is minimal.
In view of Proposition~\ref{lemtrip}(c-i)
$\T= Co(\C)$ must clearly hold whenever $\T$ is minimal. 

Conversely,  suppose that  $Co(\C) = \T$ for every section $\C$ of $\C(\T)$.
Suppose $\T$ is not minimal. Then there exists some $xy \in \T$, $x,y \in X$, such that
$\T - \{xy\}$ is a triplet cover for $T$. Thus for all $v \in \Vr$, we
have $xy$ is not contained in some $A \in S_v(\T)$. So, we can choose a section
$\C$ of $\C(\T)$ in which $xy$ is not a subset of any element of  $\C$.
But then $xy \not\in Co(\C)$, which contradicts
the assumption that $Co(\C) = \T$ for all sections $\C$ of $\C(\T)$.\\

{\em Part (c-iv):} Clearly, if $\T$ is sparse, then $\C(\T)$ contains a unique section. Conversely,
if $\C(\T)$ has a unique section, then $|\C(\T)|=|X|-2$, and so $\T$ is sparse.
\epf

\subsection{The range of $\mu(\T)$}

\begin{figure}[ht]
    \begin{center}
        \includegraphics[width=10cm]{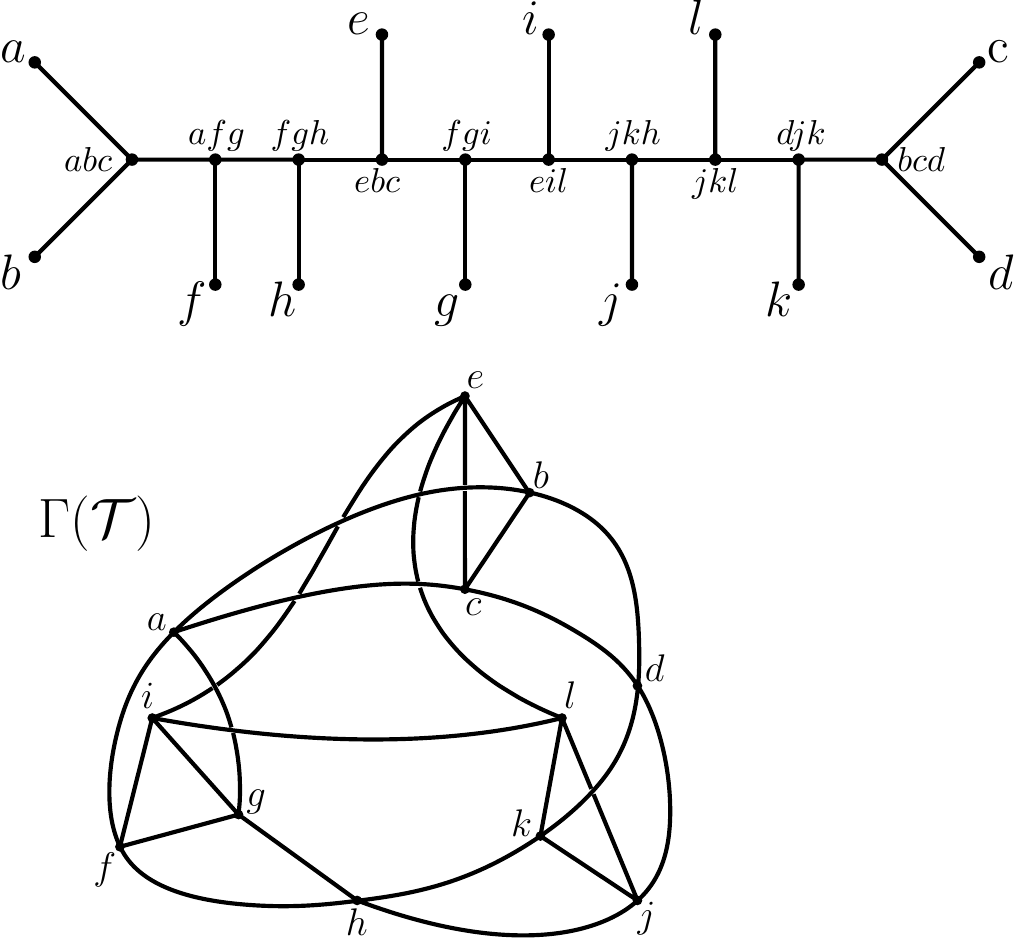}
    \end{center}
    \caption{Top:  A sparse minimal triplet cover with $\mu(\T)=4$. Bottom: The associated cover graph.}
    \label{sparse}  
\end{figure}

We saw in Lemma~\ref{lemo2} that when $\mu(\T)$ is a minimum triplet cover for $T$ then $\mu(\T)=2$.  However, for minimal triplet covers, $\mu(\T)$ can be larger. For example, there exists a sparse triplet cover $\T$ for some $T$ with $\mu(\T) = 4$
(see Fig.~\ref{sparse}; note  that in this example $\T^{-x}$ is not a triplet cover for $T-x$ for any $x \in X$).  Our main result of this section is that for a minimal triplet cover $\T$ for $T$ then $\mu(\T)$ must lie within these two extreme values.

\begin{proposition}
\label{propX5}
 If $\T$ is a minimal triplet cover for $T$, then $2 \le \mu(\T) \le 4$.
\end{proposition}

\pf
We use an argument similar to the  proof of (3) $\Rightarrow$ (2)
of \cite[Theorem 1.1]{DS09}.

Let $\C$ be a section of  $\C(\T)$, which is of Hall-type by Proposition~\ref{lemtrip}(c-ii).  
Put $n = |X|$. For $x \in X$, let $n_{\C}(x)$
be the number of triplets in $\C$ containing $x$.
If there
exists some $x\in X$ such that $n_{\C}(x)=1$, then
 $\mu_{\T}(x)=2$ by the `only if' direction of Proposition~\ref{lemtrip}(c-iiii).  
Hence $\mu(\T)=2$.
Thus we may suppose that $n_{\C}(x) > 1$ for all $x \in X$.
 Let $\Omega=\{(x,S) \in X \times \C \,:\, x \in S \}$.
Then
$$
|\Omega|= \sum_{x \in X} n_{\C}(x) \ge 2k + 3(n-k),
$$
where $k=|\{ x \in X \,:\, n_{\C}(x)=2\}|$, and, since $\C$ is a section of $C(\T)$, 
$$
|\Omega| = 3|\C| = 3(n-2).
$$
Hence $2k+3(n-k)\le 3n-6$, and so $k\ge 6$. Hence there
exists some $x \in X$ with $n_{\C}(x)=2$. Thus $\mu_{\T}(x) \le 4$
and so $\mu(\T) \le 4$ by again invoking the `only if' direction of Proposition~\ref{lemtrip}(c-iii).
\epf

\section{2-tree decompositions}

A graph $H=(V,E)$ with $|V|\geq 3$ is called a {\em 2-tree} if
there exists an ordering $v_1,v_2,\dots,v_q$ of $V$ such
that $\{v_1,v_2\} \in E$ and, for $i=3,\dots,q$, the
vertex $v_i$ has degree 2 and belongs to a
unique triangle in the subgraph induced by $H$ on
the set $\{v_1,v_2,\dots,v_i\}$ \cite[p.235]{GLM04}.
We let $\Delta(H)$ denote the set consisting of the triangles in $H$. Note
that $|E|=2|V|-3$ and $|\Delta(H)| = |V|-2$ \cite[p.227]{LM98}.

A {\em 2-tree decomposition} of a graph $H=(W,F)$ is a collection
$\{H_i=(W_i,F_i)\}_{i=1}^m$, $m \ge 1$, of subgraphs of $H$ (called {\em blocks}) such that the following hold:
\begin{itemize}
\item
 $H_i$ is a 2-tree, and its vertex set $W_i$ has size at least 3;
\item
the vertex sets of $H_1, \ldots, H_m$ cover every vertex of $H$ (i.e. $W=\bigcup_{i=1}^m W_i$);
\item
the edge sets of $H_1, \ldots, H_m$ partition the edge set of $H$
(i.e.  $F =\bigcup_{i=1}^m F_i$ and the sets $F_i$
are pairwise disjoint).  
\end{itemize}

We call the decomposition
{\em strict} if the edge set of every triangle in $H$ is a subset of some
$F_i$, $1 \le i \le m$.

Note that if $H=(V,E)$ is a 2-tree, then it can be easily shown that
any 2-tree decomposition of $H$ has one block, namely, $H$ itself.

We now state our second main theorem.

\begin{theorem}
\label{secthm}
Suppose that $\T$ is a triplet cover for $T$.
\begin{itemize}
\item[{\rm (a)}]   If $\T$ is minimal and 
$\C$ is a section of $\C(\T)$, then  $\C = \dot\cup_{i=1}^m \Delta(H_i)$ for
$\{H_i=(W_i,F_i)\}_{i=1}^m$ some 2-tree decomposition of $\Gamma(\T)$.
Moreover, $\{H_i=(W_i,F_i)\}_{i=1}^m$ is the only 2-tree decomposition of $\Gamma(\T)$
with $\C = \dot\cup_{i=1}^m \Delta(H_i)$ and, if $\T$ is sparse, then this decomposition
is a strict 2-tree decomposition. 
\item[{\rm (b)}] The following statements are equivalent:
\begin{itemize}
    \item[{\rm (i)}] $\Gamma(\T)$ has a strict 2-tree decomposition;
    \item[{\rm (ii)}]$\Gamma(\T)$ has a unique strict 2-tree decomposition;
    \item[{\rm (iii)}]$\T$ is minimal and sparse.
\end{itemize}
\item[{\rm (c)}] If $\T$ is minimal, then the following statements are equivalent:
\begin{itemize}
    \item[{\rm (i)}] $\Gamma(\T)$ has a unique 2-tree decomposition;
    \item[{\rm (ii)}] $\T$ is sparse.
\end{itemize}
\end{itemize}

\end{theorem}

\noindent{\em Proof:}
\noindent {\em Part (a):}
Suppose that $\T$ is a minimal triplet cover for $T$, and let $\C$ be a section of $\C(T)$.
We first construct a 2-tree decomposition of $\Gamma(\T)$.

Pick some element $t_1 \in \C$. We now select a sequence 
of elements from $\C$ starting with $t_1$ as follows. Suppose that some 
sequence $t_1, t_2, \dots, t_i \in \C$, $i \ge 1$,  of elements in $\C$ has been selected.
To select $t_{i+1}$ check if there exists some $t \in \C - \{t_1,t_2,\dots,t_i\} $ and
some $1 \le j \le i$ with $|t_j \cap t| = 2$. If such a $t$ exists, then put $t_{i+1}=t$,
and repeat this process for the new sequence,
otherwise stop. This process will clearly stop yielding a sequence
$t_1,t_2,\dots,t_k$ with $1 \le k \le |\C|$.

Now put $(W_1 = \bigcup_{i=1}^k t_i$, $F_1 = \bigcup_{i=1}^k Co(\{t_i\}) )$.
If $\C = \{t_1,t_2,\dots,t_k\}$ then stop. Otherwise pick some 
element in $\C - \{t_1,t_2,\dots,t_k\}$, and 
repeat the above
process for $\C - \{t_1,t_2,\dots,t_k\}$, to obtain a new pair
$(W_2,F_2)$. Then, if necessary, keep repeating this whole process until all
elements in $\C$ have been selected. This results
in a collection of pairs $\{ (W_1,F_1), \dots (W_m,F_m)\}$, with $1 \le m \le |\C|$.

Now, note that by construction $H_i=(W_i,F_i)$
is a subgraph of $\Gamma(\T)$, $|W_i| \ge 3$, and $H_i$
a 2-tree for all $1 \le i \le m$.
Moreover, since $\bigcup\C =X$ (see Proposition~\ref{lemtrip}(c-i)), it follows that $X=\bigcup_{i=1}^m W_i$,
and since $\T$ is minimal, $Co(\C)=\T$ (see Proposition~\ref{lemtrip}(c-iii)) and so $\T = \bigcup_{i=1}^m F_i$.

We now observe that if $e$ is any element of $\T$, then
by construction of the pairs $(W_i,F_i)$, there must be some $1 \le l \le m$
with $e \in F_l$ and $e \not\in F_k$ for any $k <l$. Moreover, by construction
if $t \in \C$ with $e \in Co(\{t\})$, then $Co(\{t\}) \subseteq F_l$. In particular, again by construction,
it follows that $e \not\in F_k$ for any $k > l$.
Hence, the sets $F_i$ are pairwise disjoint.
Moreover, by construction, $\C = \dot\cup_{i=1}^m \Delta(H_i)$, and so $|X|-2=|\C| = \sum_{i=1}^m|\Delta(H_i)|$
as $\C$ is a section. Hence $\{H_i\}_{i=1}^m$ is a 2-tree decomposition of $\Gamma(\T)$,
with $\C = \dot\cup_{i=1}^m \Delta(H_i)$.

To see that the uniqueness statement holds,  
suppose that $\{H'_i=(W'_i,F'_i)\}_{i=1}^q$ is any 2-tree decomposition 
of $\Gamma(\T)$ 
with $\C = \dot\cup_{i=1}^q \Delta(H_i')$. Suppose that $v_1,\dots,v_b$ is an
ordering of $W'_1$, which can be used to construct the 2-tree $H_1'$. Let $\delta$ be the triangle $ v_1,v_2,v_3$  in $H'_1$ (which
must exist as $|W_1'| \ge 3$).  Then
$\delta$ must be contained in $H_i$ for some $1 \le i \le m$, since
$\C = \dot\cup_{i=1}^m \Delta(H_i)$. 

Now, note that $H_1'$ is a subgraph of $H_i$. 
Indeed, if not then there  must exist some $4 \le l \le b$ so that
the triangle added to form a 2-tree on $v_1,\dots,v_l$ which is subgraph of $H'_1$ is not contained in $H_i$, but the
2-tree obtained from the sequence $v_1,\dots,v_{l-1}$ is in $H_i$. So the triangle containing $v_l$
which is added at stage $l$ to $H_1'$ must be contained in some $H_j$, $j \neq i$. But
this contradicts $F_i \cap F_j = \emptyset$.

Using similar reasoning, it follows that
$H_i$ is a subgraph of $H_1'$ (since we can also use the triangle $\delta$ as 
the first three elements in an ordering for
constructing $H_i$). Thus $H_1'$ is equal to $H_i$.  

Now, we can repeat this
process for $H'_2$, considering now a triangle that provides the
first three elements in an ordering for constructing $H'_2$
that is in the set $\C - \Delta(H'_1)$, and keep repeating this
whole process until finally come to considering a triangle in $H'_q$
in the set $\C - \dot\cup_{i=1}^{q-1} \Delta(H'_i)$. In this way we see that
$q=m$ and  $H'_i$ is equal to some $H_j$ for all $1\le i \le m$.
Hence the 2-tree decomposition $\{H_i\}_{i=1}^m$ is the unique
2-tree decomposition of $\Gamma(\T)$ with $\C = \dot\cup_{i=1}^m \Delta(H_i)$. 

To complete the proof of Part (a), note that if $\T$ is sparse, then $\C = \C(\T)$ by Proposition~\ref{lemtrip} (c-ii). 
It follows by Theorem~\ref{covpro}(i) and construction of $\{H_i\}_{i=1}^m$ that there is no triangle in $\Gamma(\T)$
which is not contained in some $H_i$.  Hence $\{H_i\}_{i=1}^m$
is a strict 2-tree decomposition of $\Gamma(\T)$.

\bigskip

\noindent From now on, we denote the 2-tree decomposition 
associated to a section $\C$ given in Part (a) by $\Hh_{\C}$.\\

\noindent {\em Part (b):}
The implication (ii) $\Rightarrow$ (i)  is obvious.\\
(iii) $\Rightarrow$ (i) follows from Part (a) by taking the 2-tree decomposition $\Hh_{\C(\T)}$.\\
(i) $\Rightarrow$ (iii)
Let $\{H_i=(W_i,F_i)\}_{i=1}^m$ be a strict 2-tree decomposition for  $\Gamma(\T)$.
Then $|X|-2=\sum_{i=1}^m |\Delta(H_i)|$. But every triangle in $\Gamma(\T)$
is contained in precisely one of the blocks $H_i$ since $\{H_i=(W_i,F_i)\}_{i=1}^m$ is strict.
Hence by Theorem~\ref{covpro}(i),
$$
|\C(\T)| =|\Delta(\Gamma(\T))| = \sum_{i=1}^m |\Delta(H_i)| = |X|-2,
$$
and so $\T$ is sparse.
Moreover, $\T$ is minimal. Indeed, every edge $xy$ of $\Gamma(\T)$ must be
contained in some triangle of $\Gamma(\T)$ (by definition of a 2-tree decomposition).
Hence, its removal would imply $|\C(\T-xy)| =|\Delta(\Gamma(\T-xy))| < |\Delta(\Gamma(\T))| = |X|-2$,
and so $\T-xy$ would not be a triplet cover of $T$.\\

(i) $\Rightarrow$ (ii) Suppose that
$\{H_i=(W_i,F_i)\}_{i=1}^m$ is any strict 2-tree decomposition for $\Gamma(\T)$.
Then $\C(\T) = \Delta(\Gamma(\T))=\dot\cup_{i=1}^m \Delta(H_i)$,
and since $\T$ is sparse $\C(\T)$ has a unique section (by Proposition~\ref{lemtrip}(c-iv)), namely $\C(\T)$.
Statement (ii) now follows from Part (a) and the fact that $\T$ is minimal.

\bigskip

\noindent{\em Part (c):}
(i) $\Rightarrow$ (ii) Suppose that $\T$ is not sparse. Then $\C(\T)$ must have
at least two distinct sections $\C$ and $\C'$.   Then
$\Hh_{\C} \neq \Hh_{\C'}$ are 2-tree decompositions of $\Gamma(\T)$. Statement (i)
now follows immediately.\\
(ii) $\Rightarrow$ (i): This follows from Part (b).
\epf



Note that there exists a minimal triplet cover $\T$ for some $T$, and
a 2-tree decomposition of $\Gamma(\T)$ which is not 
of the form $\Hh_{\C}$ as defined in the
proof of the last theorem for any $\C$ a section of $\C(\T)$.
Namely, take the tree and minimal triplet cover in Fig.~\ref{minimal}, add
in a new leaf $r$ to the edge adjacent to the cherry $\{p,d\}$ to get $T$, and
add the cords $rc$ and $rd$ to get $\T$. Then for the 2-tree $H_1$ 
with triangles consisting of the vertex sets
$\{p,a,d\}$, $\{a,b,d\}$, $\{a,b,c\}$, $\{b,c,q\}$, $\{a,c,w\}$, $\{d,b,z\}$ 
and the 2-tree $H_2$ consisting
of the triangle $\{r,d,c\}$,  $\{H_1,H_2\}$ is a 2-tree decomposition, but cannot arise from
a section since $abc$ and $abd$ are both in the support of the same vertex of $T$.

\subsection{Further observations}

We now show that 2-tree decompositions have some attractive
properties, which also allow us to obtain previous results on
triplet covers in a rather natural way.

\begin{proposition}
\label{y0lem}
Suppose that $H=(W,F)$ is a graph that has a 2-tree decomposition into $m\geq 1$
blocks.
Then $|F| = 2|W|-4+m$.\\
\end{proposition}
\pf
Let $\{H_i=(W_i,F_i)\}_{i=1}^m$ be a 2-tree  decomposition of $H$.
Then, since the number of triangles in a 2-tree equals the size of
its vertex set minus 2
\begin{equation}\label{one}
|W|-2 =  \sum_{i=1}^m|\Delta(H_i)| = \sum_{i=1}^m (|W_i| -2).
\end{equation}

Now, as $\{H_i\}_{i=1}^m$ is a 2-tree  decomposition,
$|F| = \dot\cup_{i=1}^m |F_i|$. Hence, as
the number of edges in a 2-tree is equal to twice the number of its vertices minus 3,
we have
\begin{equation}\label{two}
|F|  = \sum_{i=1}^m (2|W_i| -3).
\end{equation}
Using  Equations (\ref{one}) and  (\ref{two}) it immediately
follows that $|F| = 2|W|-4+m$.
\epf

Proposition~\ref{y0lem} leads directly to the following result
which shows that the size of the 2-tree decomposition 
$\Hh_{\C}$ associated to a section $\C$ is independent of the choice of $\C$.

\begin{corollary}
\label{y2cor}
Suppose that $\T$ is a minimal triplet cover for $T$, for which
$\Gamma(\T)$ has a 2-tree decomposition into $m\geq 1$ blocks. Then
$|\T| = 2|X|-4+m$. In  particular, it follows that if $\C$ is any section of
$\C(\T)$, then the 2-tree decomposition $\Hh_{\C}$
of $\Gamma(\T)$ has $|\T| -2|X|+4$ blocks.\\
\end{corollary}

We pause to mention two consequences of Corollary~\ref{y2cor}.
First, suppose that $\T$ is a triplet cover for $T$. Let $\T'$ be a minimal
triplet cover for $T$ contained in $\T$, and let $\C$ be some section of $\C(\T')$. Let $m = |\Hh_{\C}|$. Then
by Corollary~\ref{y2cor}, $|\T| \ge |\T'| = 2|X|-4 +m$. Thus, since $m \ge 1$, $|\T| \ge 2|X|-3$. (This
recovers \cite[Proposition 1]{HMS17}).

Second, suppose that $\T$ is a minimum triplet cover for $T$. Let $\C$  be a section
of $\C(\T)$, and $|\Hh_{\C}|=m$. By Corollary~\ref{y2cor},  we then have $2|X|-3 = |\T| = 2|X|-4 +m$, and so $m=1$.
Therefore, $\Gamma(\T)$ is a 2-tree (this  recovers \cite[Theorem 1]{HMS17}).
Moreover, as $\Gamma(\T)$ is a 2-tree, $|\C(\T)| = |X|-2$, and so $\T$ is sparse.

\bigskip

\noindent {\bf Remark}

\begin{itemize}

\item[]
The results above suggests the following natural question:
If a graph $H=(V,F)$ is 2-connected and $H$
has a strict 2-tree decomposition in which every pair of blocks
intersect in at most 2 vertices, then does there exist a minimal, sparse
triplet cover $\T$ for some phylogenetic $X$-tree $T$, 
with $\Gamma(\T)$ isomorphic to $H$?
Note that this can be shown to hold in case $H$ is a 2-tree.
\end{itemize}

\section{Shellings and ample patchworks}

The concept of a shellable triplet cover was introduced in \cite{DHS12}, and has proved helpful in subsequent papers.   In order to define it, one requires the notion of a quartet tree.  Suppose that $X=\{a,b,c,d\}$ and that $T$ is a phylogenetic tree for which the 
path joining $a$ and $b$ does not share a vertex with
the path joining $c$ and $d$.  In that case we say that $T$ is a 
{\em quartet tree} and denote it by writing $ab|cd$.

Given a triplet cover $\T$ of a phylogenetic $X$--tree $T$, we say that $\T$ is
{\em $T$-shellable} if either $|X|=3$ or $|X|\ge 4$ and there exists an ordering of the cords in
${X \choose 2} - \T$, say $a_1b_1, a_2b_2, \dots, a_mb_m$
such that for every $1 \le i \le m$, there exists
a pair $x_i,y_i$ of distinct elements in $X - \{a_i,b_i\}$
such that the restriction of $T$ to the
set $Y_i = \{a_i,b_i,x_i,y_i\}$ is the
quartet tree $x_ia_i | y_ib_i$, and all cords in $Y_i \choose 2$
except $a_ib_i$  are contained in $\T_i = \T \cup \{a_jb_j \,: 1 \le j \le i-1\}$.

For example, the triplet cover (indicated in terms of the median)  in Fig.~\ref{fig1} is shellable since $ae, bd, ad$ is a shelling  for it.

Although this combinatorial definition of shellability  seems somewhat involved, its motivation rests on it being a sufficient condition for recursively determining the distances between all pairs of leaves (when the edges of $T$ are assigned arbitrary positive edge lengths)  starting with just the distance values for the pairs in the triplet cover. In other words,  if a  triplet cover $\T$ of $T$ is shellable then the pairs of elements from $X$ that are not already present in $\T$ can be ordered in a sequence so  that the distance in $T$ between the leaves in each pair is uniquely determined from the distances values on pairs that are either (i) present
as an element of $\T$ or (ii) appear earlier in the sequence.

\begin{figure}[ht]
    \begin{center}
    \label{rookfig}
        \includegraphics[width=7cm]{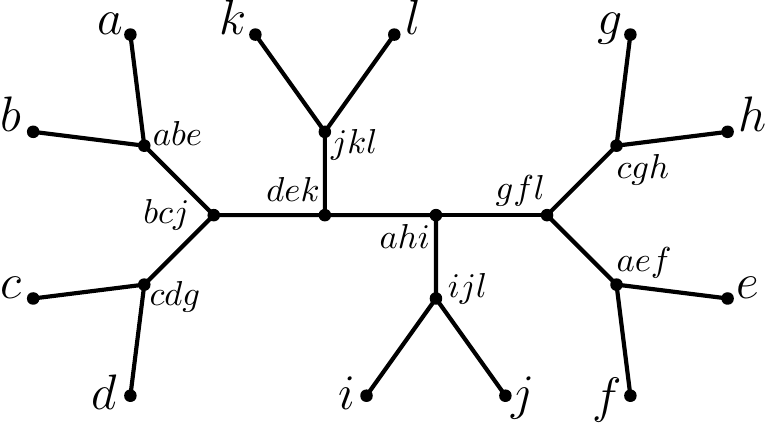} \includegraphics[width=4.5cm]{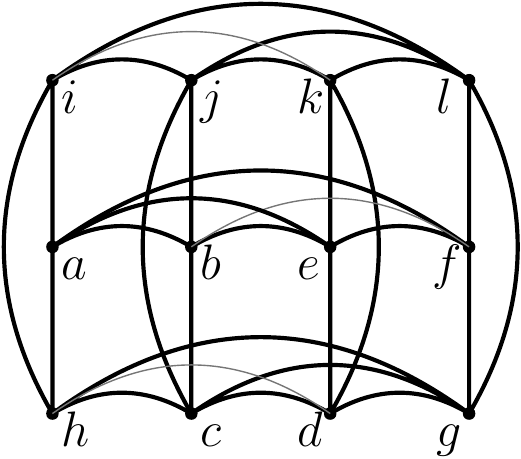}
    \end{center}
    \caption{A non-shellable triplet cover (left).  Its associated cover graph is described by the bold edges in the graph on the right. The three additional (lightly shaded) edges $ik, bf, dh$ form the initial sequence of a shelling (based on the quartets $ij|kl, ab|ef$, and $cd|gh$,  respectively), which does not extend to a full shelling. }\label{rook}
\end{figure}

Note that a triplet cover of a tree $T$ need not be $T$--shellable, even if $\T$ is sparse.  An example is shown in  Fig.~\ref{rook}.  Thus, it is of interest to better understand those triplet covers which are shellable.

The following lemma recalls some basic properties of shellability established in Proposition 4 of \cite{HMS17}.
\begin{lemma}
\label{shelem}
\mbox{}
\begin{itemize}
\item[{\rm (i)}]  Suppose that $x \in X$ and  $\T$ is a triplet cover of $T$
such that $\T^{-x}$ is a triplet cover of $T-x$.
If $\T^{-x}$ is $(T-x)$-shellable, then $\T$ is $T$-shellable.
\item[{\rm (ii)}] Suppose that $\T$, $\T'$ are triplet
covers of $T$ and that $\T' \subseteq \T$. If $\T'$ is
$T$-shellable, then so is $\T$.
\end{itemize}

\end{lemma}

Let $M$ be a finite set.   A  {\em hierarchy  on  $M$} is a collection $\Hh$ of non-empty
subsets of $M$ which satisfies the property that for all $A,B \in \Hh$, $A \cap B \in \{\emptyset,A,B\}$.
The hierarchy $\Hh$ is {\em maximal} if there is no element $H$ of $2^M - \Hh$ such that $\Hh \cup \{H\}$
is a hierarchy.

From \cite{BD01}, a collection $\Pa$ of subsets of $M$ forms a {\em patchwork}
if it satisfies the following property
$$
A,B \in \Pa \mbox{ and } A \cap B \neq \emptyset \Rightarrow A \cup B \in \Pa.
$$
A patchwork $\Pa$ on $M$  is called {\em ample} if $M$ is in $\Pa$, for
all $m \in M$, $\{m\}$ is in $\Pa$, and $\Pa$ contains a maximal hierarchy on $M$.

 If $\T$ is a triplet cover of $T$ and $\C$ is a section
of $\C(\T)$, we define $\Pa(\C)$ to be
the collection of non-empty subsets $\C'$ of $\C$ that satisfy
$|\bigcup \C'| = |\C'|+2$.
Note that $\Pa(\C)$ is a patchwork by Proposition~\ref{lemtrip}(c-ii) and \cite[Lemma 1.2]{DS09}.

Note that for the sparse minimal triplet cover $\T$ of $T$ in
Fig.~\ref{sparse}, $\T$ is $T$-shellable, but $\Pa(\C(\T))$ is not an ample patchwork.

\begin{proposition}
\label{proshel}
\mbox{}
Suppose that  $\T$ a triplet cover for $T$ and $\C$ is a section of $\C(\T)$. 
\begin{itemize}
\item[{\rm (i)}] 
If  $Co(\C)$ shellable, then so is $\T$. 
\item[{\rm (ii)}] 
If $\C',\C'' \in \Pa(\C)$ and $\C' \cap \C''= \emptyset$,
then $|\bigcup \C' \cap \bigcup\, \C''| \le 2$. Moreover, if 
$\C' \cup \C'' \in \Pa(\C)$ also holds, then $|\bigcup \C' \cap \bigcup\, \C''| = 2$.
\end{itemize}
\end{proposition}

\pf
Part (i) follows by Lemma~\ref{shelem}(ii).
For Part (ii), let $\C',\C'' \in \Pa(\C)$ be such
that  $\C' \cap \C''= \emptyset$ and let  
$Y = \bigcup \C'$ and $Z = \bigcup \C''$. Then 
since $\C$ is a section of $\C(\T)$ and 
$\C'\cup\C''\subseteq \C$ we obtain
\begin{eqnarray*}
	|Y| + |Z| -|Y \cap Z| & = & |Y \cup Z|\\
	& = & |\bigcup (\C'  \cup \C'')| \\
	& \ge & |(\C'  \cup \C'')| +2\\
	& = & | \C' | + | \C''| + 2\\
	& = & (|Y|-2)+(|Z|-2)+2 =  |Y| + |Z| - 2
\end{eqnarray*}
in view of  Lemma~\ref{lemtrip}(c--ii).
Moreover, if $\C' \cup \C'' \in \Pa(\C)$, then equality holds in the third line.
\epf

Suppose $\T$ is a triplet cover for $T$. For $A$ a non-empty subset of $X$ of size at least three,
we define $T|_A$ to be the subtree of $T$ spanned by the leaves in $A$ (suppressing degree 2 vertices). Clearly,
$T|_A$ is a phylogenetic $A$-tree. We  also
define $\T|_A$ to be the subset of $\T$ consisting of those cords $xy \in \T$
with $x,y \in A$.

\begin{lemma}
	\label{sectrip}
	Suppose $\T$  is a triplet cover for $T$ and that $\C$
	is a section of $\C(\T)$.
	\begin{itemize}
		\item[{\rm (i)}]  If $\C' \in \Pa(\C)$ then, for $Y= \bigcup\, \C'$,
		$\T|_Y$ is a triplet cover for $T|_Y$.
		\item[{\rm (ii)}]
		Suppose $\C',\C'' \in \Pa(\C)$, $A = \bigcup \C'$, $B = \bigcup \C''$,
		$A\cup B=X$ and $A \cap B = \{x,y\}$
		for some $x \neq y \in X$. If $a \in A - \{x,y\}$ and $b \in B - \{x,y\}$, then
		the quartet tree induced by $T$ on $\{a,b,x,y\}$ is either $ax|yb$ or $ay|bx$.
	\end{itemize}
\end{lemma}

\pf
{\em Part (i):}  Note that (considering $T|_Y$ as a subtree of $T$)
\begin{equation}\label{contain}
\{ v \in \Vr \,:\, v = {\rm med}_T(x,y,z) \mbox{ and } xyz \in \C'\} \subseteq \Vr(T|_Y).
\end{equation}
But $|\Vr(T|_Y)| = |Y|-2$ (since $T|_Y$ is a phylogenetic $Y$--tree), and, since
$\C$ is a section of $\C(\T)$ and $\C' \in \Pa(\C)$ we have:
$$
|\{ v \in \Vr \,:\, v = {\rm med}_T(x,y,z) \mbox{ and } xyz \in \C'\} | =|\C'| = |\bigcup \C'| -2  = |Y|-2.
$$
Therefore, equality holds in (\ref{contain}), from which Lemma~\ref{sectrip}(i) immediately follows.\\

\noindent {\em Part (ii):} Suppose that $\C',\C'' \in \Pa(\C)$, $A = \bigcup \C'$, $B = \bigcup \C''$,
$A\cup B=X$ and $A \cap B = \{x,y\}$
for some $x \neq y \in X$. If $a \in A - \{x,y\}$ and $b \in B - \{x,y\}$, then we claim that
the quartet tree induced by $T$ on $\{a,b,x,y\}$ is either $ax|yb$ or $ay|bx$.\\
\noindent {\em Proof of Claim:}
Let $v \neq x,y$ be any vertex in $T$ on the path between $x$ and $y$ in $T$.
Let $u \in V(T)$ be the vertex in $T$ that is adjacent to $v$ but not on the
path between $x$ and $y$. Consider the subtree $T_u$ of $T$ which
is the component obtained by removing the edge $\{u,v\}$ from $T$
that contains $u$. We will show that either $V(T_u)\cap X \subseteq A$ or
$V(T_u)\cap X \subseteq B$. The statement then follows immediately.

Consider $T_u$ as being a directed rooted tree with root $u$ and all
edges directed away from $u$.  For $w \in V(T_u)$, let $T_w$ denote
the directed, rooted subtree of $T_u$ with root $w$, and let $Y_w = V(T_w) \cap X$.
We show by induction on $|Y_w|$ that $Y_w \subseteq A$ or $Y_w \subseteq B$
for all $w \in V(T_u)$.

If $|Y_w|=1$, then since
$A \cup B=X$, clearly $Y_w \subseteq A$ or $Y_w \subseteq B$. Suppose that  $|Y_w| > 1$. Let $w',w''$ be the children of $w$ in $T_w$. Then, by induction, $Y_{w'}$ is a subset of $A$ or $B$ and so
is $Y_{w''}$.
Suppose without loss of generality that
$Y_{w'}$ is a subset of $A$. We need to show that
$Y_{w''}$ is also a subset of $A$.
Since $\T$ is a triplet cover of $T$, there must exist a cord
$pq\in \T$  with $p \in Y_{w'}$ and $q \in Y_{w''}$, and some $r\in X-\{p,q\}$ with
$pqr \in S_w(\T)$.   Hence, $pqr \in B$ or $pqr \in A$. But since $p \in A-\{x,y\}$ and $A \cap B = \{x,y\}$ it follows that $q \in A$.
Hence, $Y_{w''} \subseteq A$.
\epf

We can now state our third main theorem.

\begin{theorem}
	\label{patchthm}
	If $\T$ is a triplet cover of $T$ and there exists a
	section $\C$ of $\C(\T)$ such that $\Pa(\C)$ is an ample patchwork, then $\T$ is shellable.
\end{theorem}

\noindent {\em Proof:} Since $\Pa(\C)$ is an ample patchwork, there exists a maximal
hierarchy $M$ in $\Pa(\C)$. Suppose that $\C' \in \Pa(\C)$, $\C' \neq \{t\}$, some $t \in \C$,
and $\C_1, \C_2$ are
the children of $\C'$ in $M$. Let $Y= \bigcup \C'$, and $Y_i = \bigcup \C_i$, $i=1,2$,
so that in particular $Y= Y_1 \cup Y_2$. By Lemma~\ref{sectrip}(i), $\T|_Y$ is
a triplet cover for $T|_Y$, and
$\T|_{Y_i}$ is a
triplet cover of $T|_{Y_i}$, for $i=1,2$. We now show that if $\T|_{Y_i}$ is
$T|_{Y_i}$-shellable, for $i=1,2$, then $\T|_Y$ is
$T|_Y$-shellable. This will complete the proof of Theorem~\ref{patchthm}
since induction can be used in a bottom-up fashion on $M$ to see that $\T$ is shellable.

Since $\T|_{Y_i}$ is $T|_{Y_i}$-shellable, for $i=1,2$, it follows
that for each $i$ there must be some ordering of the cords in $\T|_{Y_i}$ which
satisfies the definition of $T|_{Y_i}$-shellability. Hence we
can assume that we have added all cords in $Y_i \choose 2$, $i=1,2$ into $\T|_Y$
to obtain a new set of cords $\T'$.

Now, by Proposition~\ref{proshel} (applied with $\C'=\C_1$ and $\C''=\C_2$),
$|Y_1 \cap Y_2|=2$. Let $Y_1 \cap Y_2 = \{y,z\}$, $y,z \in X$.
Note  that $yz \in \T'$ by our assumption on $\T'$. Now if
$p \in Y_1-\{y,z\}$ and $q \in Y_2-\{y,z\}$ and $pq \not\in \T'$,
then $\{yz,yp,yq,zp,zq\} \subseteq \T'$, and so by Lemma~\ref{sectrip}(ii)
applied to $T|_Y$ and the sets $Y_1,Y_2$, we can add the
cord $pq$ into $\T'$. This can be
repeated until we obtain all cords in $Y \choose 2$, in such a way
that it follows that $\T|_Y$ is $T|_Y$-shellable.\epf

\begin{corollary}
	If $\T$ is a sparse triplet cover of $T$ such that $\Pa(\C(\T))$ is ample, then
	$\T$ is shellable.
\end{corollary}


\begin{corollary}
	If $\T$ is a triplet cover for $T$ for which $\Gamma(\T)$
	has a strict 2-tree decomposition into 2 or fewer blocks then $\T$ is shellable.
\end{corollary}

\pf
First note that, by Theorem~\ref{secthm}(b) and Proposition~\ref{lemtrip}(c-iv), $\C(\T)$ is the
unique section of $\C(\T)$. Now suppose that
$\Gamma(\T)$ has one block. Then $\Gamma(\T)$  is a 2-tree.
Pick some
ordering $t_1,t_2,\dots,t_p$, $p\geq 1$, of the triangles in $\Gamma(\T)$ so that
$\Gamma(\T)$ can be constructed by adding in one triangle at a time in the
given ordering. Then
it is straightforward to see that the set
$$
\{ \{t_1\},\dots,\{t_p\}, \{t_1,t_2\}, \{t_1,t_2,t_3\},\dots,\{t_1,t_2,\dots,t_p\} \}
$$
is a maximal hierarchy in $\Pa(\C(\T))$. Now apply Theorem~\ref{patchthm}.

If the 2-tree decomposition of $\Gamma(\T)$ has two blocks,
let $t_1,t_2,\dots,t_p$ and $t'_1,t'_2,\dots,t'_q$ be orderings
of the triangles in each of the blocks  so that each block
can be constructed by adding in one triangle at a time in the
given ordering. Put $\C'=\{ \{t_1\},\dots,\{t_p\}, \{t_1,t_2\}, \{t_1,t_2,t_3\},\dots,\{t_1,t_2,\dots,t_p\}\} $, $p\geq 1$,
and $\C''= \{ \{t'_1\},\dots,\{t'_q\}, \{t'_1,t'_2\}, \{t'_1,t'_2,t'_3\},\dots,\{t'_1,t'_2,\dots,t'_q\} \}$, $q\geq 1$.
Note that $\{t_1,t_2,\dots,t_p\} \cup \{t'_1,t'_2,\dots,t'_q\} \in \Pa(\C(\T))$ by
Theorem~\ref{covpro}(ii) and  Proposition~\ref{proshel}(ii).
%
Now it is straightforward to see that the set
$$
\C'\cup\C''\cup\{\{ t_1,t_2,\dots,t_p, t_1',t_2', \dots,t_q'\} \}
$$
is a maximal hierarchy in $\Pa(\C(\T))$. Now apply Theorem~\ref{patchthm} again.
\epf

Note that the case for one block in the last corollary was also shown to hold in \cite{HMS17}.

\section{\bf Acknowledgements}
KTH and VM thank the London Mathematical Society and the
Biomathematcs Research Center, University of Canterbury,
for their support. MS thanks the (former) Allan Wilson Centre for funding support.

\section*{References}

\bibliography{triplet-cover}

\end{document}